\documentclass[a4paper]{amsart}
\usepackage{amssymb,amsthm,amsmath}
\usepackage{times,mathptmx} 
\address{Department of Mathematics, Faculty of Education and Human Sciences, 
Yokohama National University, 240-8501 Japan.}
\email{takase@ynu.ac.jp}
\author{Masamichi Takase}
\title{The Hopf invariant of a Haefliger knot}
\date{}
\newcommand{\Zset}{\mathbf{Z}}
\newcommand{\Rset}{\mathbf{R}}
\newcommand{\Cset}{\mathbf{C}}
\newcommand{\Int}{\mathop{\mathrm{Int}}\nolimits}
\newcommand{\Imm}{\mathop{\mathrm{Imm}}\nolimits}
\newcommand{\co}{\colon\thinspace}

\newcommand{\SubS}{\S}
\newtheorem{thm}{Theorem}[section]
\newtheorem{cor}[thm]{Corollary}

\newtheorem{prop}[thm]{Proposition}
\theoremstyle{definition}
\newtheorem{defn}[thm]{Definition}
\newtheorem{rem}[thm]{Remark}
\def\ack{\section*{Acknowledgements}}
\begin{document}\sloppy
\maketitle
\begin{abstract}
We show that 
Haefliger's differentiable $(6,3)$-knot bounds, in $6$-space, 
a $4$-manifold (a Seifert surface) of 
arbitrarily prescribed signature. 
This implies, according to our previous paper, 
that the Seifert surface has been prolonged 
in a prescribed direction near its boundary. 
This aspect enables us to understand a resemblance between 
Ekholm-Sz\H{u}cs' formula for the Smale invariant and 
Boech\'at-Haefliger's formula for Haefliger knots. 
As a consequence, we show that 
an immersion of the $3$-sphere in $5$-space can be regularly homotoped 
to the projection of an embedding in $6$-space 
if and only if its Smale invariant is even. 
We also correct a sign error in our previous paper: 
``A geometric formula for Haefliger knots'' [Topology 43 (2004) 1425-1447]. 
\end{abstract}
\section{Introduction}
Haefliger \cite{hae1,hae2} discovered 
differentiably knotted spheres in codimension greater than three. 
In a particular case, he proved the isomorphism 
\[\Omega\co C^{2k+1}_{4k-1}\to\Zset,\]
where $C^{2k+1}_{4k-1}$ denotes the group of 
differentiable isotopy classes of embeddings 
of $S^{4k-1}$ in $S^{6k}$ ($k\ge1$). 

We deal with the case where $k=1$. 
We know that any such Haefliger knot $F\co S^3\hookrightarrow S^6$ as above 
has a Seifert surface, that is, an embedding $\widetilde{F}\co V^4\hookrightarrow S^6$ of 
a compact oriented $4$-manifold whose restriction to the boundary $\partial{V^4}(=S^3)$ 
coincides with $F$ 
(see \cite{boechat,b-h},\cite[p.95]{g-m} and \cite[Corollary~6.2]{takase2}). 
Furthermore, we can consider, associated to such a Seifert surface, 
the Hopf invariant $H_{\widetilde{F}}\in\Zset$ of the map 
from $S^3$ into $S^6\smallsetminus F(S^3)\simeq S^2$, 
determined by the outward normal field of $F(S^3)\subset\widetilde{F}(V^4)$. 
Then, by \cite[Corollary~6.5]{takase2}, 
the Haefliger invariant $\Omega(F)$ of $F$ is represented as 
\[
\Omega(F)=-\frac{1}{8}(\sigma(V^4)+H_{\widetilde{F}}), 
\]
where $\sigma(V^4)$ denotes the signature of $V^4$. 

Regarding this formula, 
we pose a natural question: 
\textit{for a given Haefliger knot, which integer 
can be realised as the Hopf invariant of a Seifert surface?}
In Corollaries~\ref{cor:hopf} and \ref{cor:signature}, 
we give a complete answer to it; 
namely, we show that given a Haefliger knot and given an integer, 
we can choose a Seifert surface 
so that its Hopf invariant is equal to 
(or its signature is equal to) the given integer. 

In \SubS\ref{subsect:erratum}, we will correct a sign error of 
our previous paper \cite[Lemma~5.2]{takase2}. 
With this correction, 
the sign should be changed 
in each term involving the square of the normal Euler class, 
appearing in \cite[Theorem~5.1, Corrollaries~6.2, 6.3(a) and 6.5]{takase2}. 
For details, see \SubS\ref{subsect:erratum}. 

As a consequence, 
with the help of Ekholm and Sz\H{u}cs' formula for 
the Smale invariant $\omega\co\Imm[S^3,\Rset^5]\to\Zset$ 
(denoting by $\Imm[S^3,\Rset^5]$ the group of regular homotopy classes 
of immersions of $S^3$ in $\Rset^5$), 
we show that 
for an embedding $F\co S^3\hookrightarrow\Rset^6$ 
and an immersion $f\co S^3\looparrowright\Rset^5$ with even Smale invariant, 
we can isotope $F$ so that its composition 
with the projection $\Rset^6\to\Rset^5$ 
becomes an immersion regularly homotopic to $f$ (Theorem~\ref{thm:proj}). 
In particular, any immersion $S^3\looparrowright\Rset^5$ with even Smale invariant 
can be regularly homotoped to the projection of the unknot in $C^3_3$. 
In Corollary~\ref{cor:proj}, we show that 
an immersion $f\co S^3\looparrowright\Rset^5$ 
can be regularly homotoped to the projection of an embedding 
$S^3\hookrightarrow\Rset^6$ if and only if its Smale invariant 
$\omega(f)$ is even. 

In \S\ref{sect:punc}, we discuss aspects of embeddings of 
punctured $4$-manifolds in $S^6$. 

We work in the $C^\infty$-differentiable category; 
all manifolds and mappings 
are supposed to be differentiable of class $C^\infty$, 
unless otherwise explicitly stated. 
Throughout this paper, 
all $C^\infty$-mappings shall be in 
the so called \textit{nice dimensions} (see \cite{mather}), 
where the two notions \textit{stable map} and \textit{generic map} are equivalent. 
We will use the term \textit{generic} according to \cite[\SubS 2.1]{e-s}. 
For a map $f\co X\to Y$ from a manifold $X$ with non-empty boundary 
we often denote its restriction to the boundary by 
$\partial{f} (:=f|_{\partial{X}}\co\partial{X}\to Y$). 

We will suppose the spheres are oriented. 
If $M$ is an oriented manifold with non-empty boundary, 
then for the \textit{induced orientation of\/ $\partial{M}$} 
we adopt the \textit{outward vector first} convention: 
we say an ordered basis of 
$T_p(\partial{M})$ ($p\in\partial{M}$) is positively oriented 
if an outward vector followed by the basis
is a positively oriented basis of $T_pM$. 
For a closed $n$-dimensional manifold $M$ we denote its punctured manifold 
by $M_\circ$; i.e., $M_\circ:= M\smallsetminus\Int{D^n}$.

We use the symbol `$\approx$' for a group isomorphism and 
`$\cong$' for a diffeomorphism between manifolds; 
the symbol `$\simeq$' means 
a homotopy equivalence between two topological spaces. 
The homology and cohomology groups are supposed to be 
with integer coefficients unless otherwise explicitly noted. 

\section{Seifert surfaces for Haefliger knots and their Hopf invariants}
Haefliger 
has proved in \cite{hae1} and \cite[\SubS5.16]{hae2} that 
the group $C^{2k+1}_{4k-1}$ of isotopy classes of embeddings 
$S^{4k-1}\hookrightarrow S^{6k}$ is isomorphic to the integers $\Zset$
($k\ge1$). 
Since he has also given an explicit construction of 
an embedding representing a generator of $C^{2k+1}_{4k-1}$, 
with respect to his generator 
we have the identification of $C^{2k+1}_{4k-1}$ with $\Zset$, 
which we call the Haefliger invariant 
\[
\Omega\co C^{2k+1}_{4k-1}\to\Zset.
\] 

For a Haefliger knot, 
we can consider its \textit{Seifert surface} 
analogously as in the usual knot theory in codimension two. 

\begin{defn}\label{def:seifert}
For an embedding $F\co S^{4k-1}\hookrightarrow S^{6k}$, 
\textit{a Seifert surface for $F$} is 
an embedding $\widetilde{F}\co V^{4k}\hookrightarrow S^{6k}$ 
of a compact connected oriented $4k$-manifold with boundary $\partial{V^{4k}}=S^{4k-1}$, 
whose restriction 
$\partial{\widetilde{F}}\co S^{4k-1}\hookrightarrow S^{6k}$
to the boundary coincides with $F$. 
\end{defn}

We will discuss the existence of such a Seifert surface
in \SubS\ref{subsect:constructions}. 
Associated to a Seifert surface, 
we can consider the Hopf invariant as in \cite{takase2}. 

\begin{defn}\label{def:hopf}
Let $\widetilde{F}\co V^{4k}\hookrightarrow S^{6k}$ 
be a Seifert surface for an embedding 
$F\co S^{4k-1}\hookrightarrow S^{6k}$. 
Then, the outward normal field of 
$F(S^{4k-1})\subset\widetilde{F}(V^{4k})$ 
determines the map 
$\nu_{\widetilde{F}}\co S^{4k-1}\to S^{6k}\smallsetminus F(S^{4k-1})$. 
We consider its homotopy class 
$[\nu_{\widetilde{F}}]$ 
to be lying in 
$\pi_{4k-1}(S^{6k}\smallsetminus F(S^{4k-1}))=\pi_{4k-1}(S^{2k})$, 
via the homotopy equivalence 
$p\circ\psi\co S^{6k}\smallsetminus F(S^{4k-1})\stackrel{\cong}{\to}D^{4k}\times S^{2k}\stackrel{\simeq}{\to}S^{2k}$, 
where $\psi$ is an orientation-preserving diffeomorphism 
(see \cite[Theorem~5.2]{smale}) and $p$ is the projection. 
Thus, we define \textit{the Hopf invariant} 
$H_{\widetilde{F}}$ \textit{for} $\widetilde{F}$ 
to be the Hopf invariant of $[\nu_{\widetilde{F}}]$. 
\end{defn}

\subsection{A formula for the Haefliger invariant --- a correction 
to the paper ``A geometric formula for Haefliger knots'' 
[Topology 43 (2004) 1425-1447]}\label{subsect:erratum}
In our previous paper \cite{takase2}, 
we related the Hopf invariant for $\widetilde{F}$ 
to the normal Euler class of $\widetilde{F}$; 
however, there was an error about a sign. 
This subsection is devoted to correcting the sign error 
and reviewing a formula for the Haefliger invariant 
with the corrected sign (see also \cite{boechat,b-h,g-m}). 

The error occurs in the last sentence of the proof of Lemma~5.2 
in \cite[p.1443]{takase2}. 
It reads ``\textit{Thus, we have 
$H_{\widetilde{F}}=-[\widetilde{F}(\Sigma^{2k})].$}'' but 
should be 
``\textit{Thus, the desired functional product is computed to be equal 
to $-[\widetilde{F}(\Sigma^{2k})]$ and hence the Hopf invariant 
$H_{\widetilde{F}}$ is equal to $[\widetilde{F}(\Sigma^{2k})]\in H_{2k}(X)=\Zset$}.'' 
The reason is that the Hopf invariant is defined, 
in \cite[\SubS2.2]{takase2}, 
to be minus the functional cup product. 
Note that the similar argument can be found 
in the last part of the proof of Proposition~3.5 in \cite[p.1436]{takase2}. 

With this correction, each term involving the square of the normal Euler class, 
appearing in \cite[Theorem~5.1, Corrollaries~6.2, 6.3(a) and 6.5]{takase2}, 
should change its sign, as follows. 

First, Theorem~5.1 in \cite[p.1442]{takase2} should be written as: 

\begin{thm}[{\cite[Theorem~5.1]{takase2}}]\label{thm:hopf4k}
Let $\widehat{V}^{4k}$ be a closed oriented $4k$-manifold and put 
$V^{4k}:=\widehat{V}^{4k}\setminus\Int D^{4k}$. 
Let $\widetilde{F}\co V^{4k}\hookrightarrow S^{6k}$ be an embedding 
with normal Euler class 
$e_{\widetilde{F}}\in H^{2k}(V^{4k})=H^{2k}(\widehat{V}^{4k})$. 
Then, the Hopf invariant $H_{\widetilde{F}}$ \textup{(}along the boundary\textup{)} 
is equal to 
$-e_{\widetilde{F}}\smile e_{\widetilde{F}}\in H^{4k}(\widehat{V}^{4k})=\Zset$. 
\end{thm}

Putting $M^{4k}:=S^{2k}\times S^{2k}\setminus\Int{D^{4k}}$, 
Corollary~6.2 in \cite[p.1444]{takase2} should be as: 

\begin{cor}[{\cite[Corollary~6.2]{takase2}}]\label{cor:s2kxs2k}
Let $\widetilde{E}_{a,b}\co M^{4k}\hookrightarrow S^{6k}$ be an embedding 
with normal Euler class 
$(2a,2b)\in H^{2k}(M^{4k})=H^{2k}(S^{2k}\times S^{2k})\approx\Zset\oplus\Zset$. 
Then $E_{a,b}:=\widetilde{E}_{a,b}|_{\partial M^{4k}}\co S^{4k-1}\hookrightarrow S^{6k}$
represents $ab\in C^{2k+1}_{4k-1}=\Zset$.
In particular, 
$E_{\pm1,\pm1}\co S^{4k-1}\hookrightarrow S^{6k}$
represents the generator of $C^{2k+1}_{4k-1}$.
\end{cor}

Furthermore, Corollary~6.3(a) in \cite[p.1444]{takase2} should be as: 

\begin{cor}[{\cite[Corollary~6.3(a)]{takase2}}]\label{cor:main}
\textup{(a)}
For an arbitrary embedding 
$F\co S^{4k-1}\hookrightarrow S^{6k}$, 
there exists an embedding 
$\widetilde{F}\co V^{4k}\hookrightarrow S^{6k}$
of a compact oriented $4k$-manifold $V^{4k}$ with $\partial V^{4k}=S^{4k-1}$ 
such that $\widetilde{F}|_{\partial V^{4k}}=F$. 
Furthermore, 
\begin{eqnarray*}
\Omega(F)
&=&-\frac{1}{24}(-\bar{p}_k[\widehat{V}^{4k}]+3H_{\widetilde{F}})\\
&=&-\frac{1}{24}(-\bar{p}_k[\widehat{V}^{4k}]-3e_{\widetilde{F}}\smile e_{\widetilde{F}})
\end{eqnarray*} 
gives the isomorphism $\Omega\co C^{2k+1}_{4k-1}\longrightarrow\Zset$, 
where $e_{\widetilde{F}}\in H^{2k}(V^{4k})=H^{2k}(\widehat{V}^{4k})$ 
is the normal Euler class for $\widetilde{F}$ 
and $e_{\widetilde{F}}\smile e_{\widetilde{F}}\in H^{4k}(\widehat{V}^{4k})=\Zset$ 
is the cup product. 
\end{cor}

Corollary~6.5 in \cite[p.1445]{takase2} in the case $C^3_3$ ($k=1$)
should be as: 

\begin{cor}[{\cite[Corollary~6.5]{takase2}, see also \cite[p.95]{g-m}}]\label{cor:3dim}
Every embedding $F\co S^3\hookrightarrow S^6$ extends to 
an embedding $\widetilde{F}\co V^4\hookrightarrow S^6$ 
of a compact oriented $4$-manifold $V^4$, and 
\begin{eqnarray*}
\Omega(F)
&=&-\frac{1}{8}(\sigma(V^4)+H_{\widetilde{F}})\\
&=&-\frac{1}{8}(\sigma(V^4)-e_{\widetilde{F}}\smile e_{\widetilde{F}})
\end{eqnarray*} 
gives the isomorphism $\Omega\co C^3_3\longrightarrow\Zset.$ 
\end{cor}

\subsection{Constructions of Seifert surfaces}\label{subsect:constructions} 
From now on, we will confine ourselves to the case $C^3_3$ ($k=1$). 
In this subsection, 
we review some concrete constructions of Seifert surfaces for 
Haefliger knots, 
which have been used in \cite[p.95]{g-m} and \cite[Corollary~6.2]{takase2}. 

First, we write down the particular case of Corollary~\ref{cor:s2kxs2k}
in $k=1$. 

\begin{cor}\label{cor:s2xs2}
Let $\widetilde{E}_{a,b}\co (S^{2}\times S^{2})_\circ\hookrightarrow S^{6}$ be an embedding 
with normal Euler class 
$(2a,2b)\in H^{2}((S^{2}\times S^{2})_\circ)=H^{2}(S^{2}\times S^{2})\approx\Zset\oplus\Zset$. 
Then $E_{a,b}:=\widetilde{E}_{a,b}|_{\partial (S^{2}\times S^{2})_\circ}\co S^{3}\hookrightarrow S^{6}$
represents $ab\in C^{3}_{3}=\Zset$.
\end{cor}

We see that for any (isotopy class of) Haefliger knot 
we can construct its Seifert surface by choosing a suitable pair 
$(a,b)$ of integers. 
In particular, if we put $(a,b)=(1,1)$ in Corollary~\ref{cor:s2xs2}, 
the embedding $\widetilde{E}_{1,1}\co(S^2\times S^2)_\circ\hookrightarrow S^6$ 
is a Seifert surface for the generator of $C^3_3=\Zset$, 
of signature zero and with Hopf invariant $-8$ (by Theorem~\ref{thm:hopf4k}). 

A clue given for \cite[Exercise~3(ii) on p.95]{g-m} 
elicits another construction. 

\begin{prop}[{see \cite[p.95]{g-m}}]\label{prop:cp2}
The punctured complex projective plane ${\Cset{P}^2}_{\circ}$ 
can be embedded in $S^6$ 
with normal Euler class $2k+1\in H^2({\Cset{P}^2}_{\circ})\approx\Zset$ 
$(k\in\Zset)$. 
If we denote by 
$\widetilde{E}_{k}\co{\Cset{P}^2}_{\circ}\hookrightarrow S^6$ 
such an embedding with normal Euler class $2k+1$, 
then its Hopf invariant is equal to $-4k^2-4k-1$ and 
$\partial\widetilde{E}_{k}\co S^3\hookrightarrow S^6$ 
represents $k(k+1)/2\in C^3_3=\Zset$. 
\end{prop}

\begin{proof}
Consider the $2$-dimensional disc bundle $\eta$ over ${\Cset{P}^2}_{\circ}$ 
with Euler class $2k+1\in H^2({\Cset{P}^2}_{\circ};\Zset)=\Zset$. 
Since the second Stiefel-Whitney class 
$w_2(\eta)\in H^2({\Cset{P}^2}_{\circ};\Zset_2)=\Zset_2$ 
of $\eta$ is congruent to $2k+1$ modulo $2$, 
we have $w_2(T{\Cset{P}^2}_{\circ}\oplus\eta)
=w_2(T{\Cset{P}^2}_{\circ})\smile1+1\smile1
=0\in H^2({\Cset{P}^2}_{\circ};\Zset_2)$. 
Therefore, together with $\pi_3(SO(6))=0$, 
we see that $T{\Cset{P}^2}_{\circ}\oplus\eta$ is trivial. 
Then, Hirsch's h-principle \cite{hirsch} implies that 
${\Cset{P}^2}_{\circ}$ can be immersed in $S^6$ 
with normal bundle isomorphic to $\eta$. 
Let $\widetilde{F}\co{\Cset{P}^2}_{\circ}\looparrowright S^6$ 
be such an immersion with normal Euler class $2k+1$. 
(Actually the above argument implies that 
the normal Euler class of such an immersion 
ought to be of the form 
$2k+1\in H^2({\Cset{P}^2}_{\circ};\Zset)=\Zset$). 

Since ${\Cset{P}^2}_{\circ}$ is the total space of 
a $2$-dimensional disc bundle over $\Cset{P}^1\cong S^2$, 
the immersion $\widetilde{F}$ is regularly homotopic to an embedding. 
Since a regular homotopy does not change the normal bundle, 
the first part of the proposition has been proved. 

Let 
$\widetilde{E}_{k}\co{\Cset{P}^2}_{\circ}\hookrightarrow S^6$ 
be an embedding with normal Euler class $2k+1$. 
Then, by Theorem~\ref{thm:hopf4k}, 
its Hopf invariant $H_{\widetilde{F}}$ is equal to 
$-(2k+1)\smile(2k+1)=-4k^2-4k-1\in H^4(\Cset{P}^2;\Zset)=\Zset$, where 
we consider $2k+1\in H^2({\Cset{P}^2}_{\circ};\Zset)$ to be in $H^2(\Cset{P}^2;\Zset)$. 
Thus, by Corollary~\ref{cor:3dim}, we have 
$\Omega(\partial\widetilde{E}_{k})=-(1-4k^2-4k-1)/8=k(k+1)/2$. 
This completes the proof. 
\end{proof}

\begin{rem}\label{rem:cp2bar}
In Proposition~\ref{prop:cp2}, 
if we use the complex projective plane $\overline{\Cset{P}^2}$ 
with the reversed orientation, 
we have a Seifert surface, 
of signature $-1$ with the Hopf invariant $4k^2+4k+1$, 
for an embedding representing $-k(k+1)/2\in C^3_3=\Zset$. 
\end{rem}

\begin{rem}\label{rem:any}
Corollary~\ref{cor:s2xs2} implies that 
the punctured $S^2\times S^2$ can be a Seifert surface for 
any (isotopy class of) Haefliger knot.  
On the other hand, Proposition~\ref{prop:cp2} shows that 
the punctured complex projective plane ${\Cset{P}^2}_\circ$ 
(resp. $\overline{\Cset{P}^2}_\circ$) 
can be a Seifert surface 
only for non-negative classes 
(resp. for non-positive classes) of $C^3_3=\Zset$. 
\end{rem}

\subsection{Seifert surfaces with given Hopf invariants}\label{subsect:seifert}
Regarding the formula 
in Corollary~\ref{cor:3dim}, 
it is natural to ask 
which integer can be realised 
as the Hopf invariant (or as the signature) of a Seifert surface 
for a given Haefliger knot. 
We will give a complete answer to it 
in Corollaries~\ref{cor:hopf} and \ref{cor:signature}. 

\begin{prop}\label{prop:cp2s2s2}
For the standard inclusion $S^3\subset S^6$, 
there exists a Seifert surface 
whose Hopf invariant is equal to $1$.  
\end{prop}

\begin{proof}
If we put $k=0$ in Proposition~\ref{prop:cp2}, 
then we have the embedding with Hopf invariant $-1$, denoted by 
$\widetilde{P}\co\Cset{P}^2_{\circ}\hookrightarrow S^6$, 
whose restriction to the boundary 
represents $-(1-1)/8=0\in C^3_3=\Zset$ and hence is isotopic to the standard inclusion. 
Since a differentiable isotopy implies a differentiable ambient isotopy, 
the embedding $\widetilde{P}$ 
followed by a suitable diffeomorphism on $S^6$ 
gives a desired Seifert surface for the standard inclusion. 
\end{proof}

\begin{rem}\label{rem:cp2s2s2}
The similar argument for $\overline{\Cset{P}^2}_{\circ}$ 
(putting $k=0$ in Remark~\ref{rem:cp2bar}) provides 
another Seifert surface 
$\widetilde{Q}\co\overline{\Cset{P}^2}_{\circ}\hookrightarrow S^6$, 
of signature $-1$ and with Hopf invariant $1$, 
for the standard inclusion $S^3\subset S^6$. 
\end{rem}

Proposition~\ref{prop:cp2s2s2} and Remark~\ref{rem:cp2s2s2} 
imply that 
for a given Haefliger knot 
we can change its Seifert surface 
so as to have a desired Hopf invariant 
(in return for changing the signature), 
by taking the boundary connected sum with $\widetilde{P}$ 
or with $\widetilde{Q}$. 
Namely we have: 

\begin{cor}\label{cor:hopf}
For an arbitrary embedding $F\co S^3\hookrightarrow S^6$ and 
an arbitrary integer $k\in\Zset$, 
there exists a Seifert surface for $F$ 
whose Hopf invariant is equal to $k$. 
\end{cor}

\begin{proof}
Take an arbitrary Seifert surface for $F$. 
Then, 
by taking the boundary connected sum with $\widetilde{P}$ 
(resp.\ with $\widetilde{Q}$), 
we obtain a new Seifert surface with signature which differs 
by $-1$ (resp.\ by $1$) 
to the initial one, 
without changing the isotopy class of 
the embedding $S^3\hookrightarrow S^6$ on the boundary. 
Repeat this procedure until we have a Seifert surface with 
the Hopf invariant equal to $k$. 
Thus, composing it with a suitable ambient isotopy if necessary, 
we obtain a desired Seifert surface. 
\end{proof}

In view of the formula for the Haefliger invariant 
(Corollary~\ref{cor:3dim}), 
a very similar argument provides the following. 

\begin{cor}\label{cor:signature}
For an arbitrary embedding $F\co S^3\hookrightarrow S^6$ and 
an arbitrary integer $k\in\Zset$, 
there exists a Seifert surface for $F$ 
with signature $k$. 
\end{cor}

Recall that the normal bundle of 
an embedding $F\co S^3\hookrightarrow S^6$ 
is trivial \cite{kervaire} and that 
the homotopy classes of its normal framings are 
classified by the elements of $\pi_3(SO(3))\approx\Zset$. 
Furthermore, 
the homomorphism $\pi_3(SO(3))\to\pi_3(S^2)$ 
induced by the natural map $SO(3)\to SO(3)/SO(2)=S^2$ 
and the (usual) Hopf invariant $\pi_3(S^2)\to\Zset$ 
are both isomorphisms. 
Therefore, the following is just an interpretation of 
Corollary~\ref{cor:hopf}. 

\begin{cor}\label{cor:normal}
For an arbitrary embedding $F\co S^3\hookrightarrow S^6$ and 
an arbitrary normal vector field to $F(S^3)$ in $S^6$, 
there exists a Seifert surface $\widetilde{F}\co V^4\hookrightarrow S^6$ for which 
the normal field along the boundary $F(S^3)\subset\widetilde{F}(V^4)$ is 
in accordance with the given normal vector field. 
\end{cor}

\section{The projection of a Haefliger knot} 
The purpose of this section is to prove the following two 
Theorems~\ref{thm:proj} and \ref{thm:proj2}. 
Let $p\co\Rset^6\to\Rset^5$ be the projection, 
defined by dropping off the last coordinate: 
$(x_0,\ldots,x_5,x_6)\mapsto (x_0,\ldots,x_5)$.

\begin{thm}\label{thm:proj}
Let $f\co S^3\looparrowright\Rset^5$ be an immersion 
with even Smale invariant 
and $F\co S^3\hookrightarrow\Rset^6$ an embedding. 
Then, 
there exist an immersion $f'\co S^3\looparrowright\Rset^5$ 
regularly homotopic to $f$
and an embedding $F'\co S^3\hookrightarrow\Rset^6$ 
isotopic to $F$ such that $f'=p\circ F'$. 
\end{thm}

\begin{thm}\label{thm:proj2}
If the projection $p\circ F$ of 
an embedding $F\co S^3\hookrightarrow\Rset^6$ is an immersion, 
then its Smale invariant $\omega(p\circ F)$ is even. 
\end{thm}

As an easy corollary of Theorem~\ref{thm:proj}, we have: 
\begin{cor}\label{cor:unknot}
Any immersion $f\co S^3\looparrowright\Rset^5$ 
with even Smale invariant can be regularly homotoped to 
the projection of the unknot in $C^3_3$. 
\end{cor}

Combining Theorem~\ref{thm:proj} and Theorem~\ref{thm:proj2}, we have the following. 
\begin{cor}\label{cor:proj}
An immersion $f\co S^3\looparrowright\Rset^5$ can be 
regularly homotoped to the projection of 
an embedding $S^3\hookrightarrow\Rset^6$ 
if and only if its Smale invariant $\omega(f)$ is even. 
\end{cor}

\begin{rem}
It has been shown in \cite{s-s} 
that if a closed $n$-manifold ($n\ge3$) is non-orientable 
or odd-dimensional 
then there exists an immersion with normal crossings 
$M^n\looparrowright\Rset^{2n-1}$ 
which can never be lifted to an embedding into $\Rset^{2n}$. 
\end{rem}

The proofs of Theorems~\ref{thm:proj} and \ref{thm:proj2} will depend on 
Ekholm and Sz\H{u}cs' geometric formula \cite{e-s}, 
which enables us to read off 
the Smale invariant of, hence the regular homotopy class of, 
an immersion of $S^{4k-1}$ in $\Rset^{4k+1}$ ($k\ge 1$) 
in terms of its \textit{singular Seifert surface}. 
We will first review their formula 
(only in our case, i.e., when $k=1$) in \SubS\ref{subsect:e-s} 
and then give the proofs of Theorems~\ref{thm:proj} and \ref{thm:proj2}  
in \SubS\ref{subsect:proof}. 

\subsection{Ekholm and Sz\H{u}cs' formula}\label{subsect:e-s}
We need to recall some definitions from \cite[\SubS 2.4]{e-s}. 
Let $g\co M^4\to\Rset^5$ be a generic map of a compact oriented $4$-manifold. 
Then, the set $\widetilde\Sigma(g)$ of its singular points is 
a $2$-dimensional submanifold of $M^4$, 
in which finitely many isolated \textit{cusp points} lie. 
Note that each cusp point has 
the naturally induced orientation \cite[Appendix A]{szucs2}. 

\begin{defn}[{see \cite[Definition~2.10 and Remark~6.2]{e-s}}]\label{defn:euler}
Let $g\co M^4\to\Rset^5$ be a generic map of an oriented $4$-manifold. 
Then let $\xi(g)$ denote the $2$-dimensional vector bundle over $\widetilde\Sigma(g)$, 
whose fibre over $p\in\widetilde\Sigma(g)$ is $\xi(g)_p=T_{g(p)}\Rset^5/dg(T_pM^4)$ and 
$e[\xi(g)]$ denote its Euler number (in the sense of  \cite[Definition~2.10]{e-s}). 
\end{defn}

\begin{rem}[{\cite[Remark~6.4]{e-s}}]\label{rem:cusp} 
When a generic map $g\co M^4\to\Rset^5$ of 
an oriented $4$-manifold $M^4$ is 
the composition of an immersion $G\co M^4\looparrowright\Rset^6$ 
with the projection $\Rset^6\to\Rset^5$, 
we have a more convenient description of 
the Euler number $e[\xi(g)]$ in Definition~\ref{defn:euler}.
In such a situation, 
the homology class of $\widetilde\Sigma(g)$ in $M^4$ is dual to 
the normal Euler class of $G$, 
and the bundle $\xi(g)$ is isomorphic to the normal bundle of $G$ 
restricted to $\widetilde\Sigma(g)$, 
whose Euler class evaluated by the orientation class $[\widetilde\Sigma(g)]$ 
precisely equals $e[\xi(g)]$. 
\end{rem}

Ekholm and Sz\H{u}cs \cite[Theorem~1.1 and Remark~3.1]{e-s} have given 
a formula for the Smale invariant \cite{smale} 
\[\omega\co\Imm[S^3,\Rset^5]\to\Zset,\] 
which gives a group isomorphism between 
the group $\Imm[S^3,\Rset^5]$ of regular homotopy classes 
of immersions of $S^3$ in $\Rset^5$ and the integers $\Zset$. 

\begin{thm}[{\cite[Theorem~1.1(a) and Remark~3.1]{e-s}}]\label{thm:e-s} 
Let $f\co  S^3 \looparrowright \Rset^5$ be an immersion and 
$\widetilde{f}\co V^4\to\Rset^5$ be 
a \textit{singular Seifert surface} for $f$, that is, 
a generic map with $\partial{\widetilde{f}}=f$
which has no singularity near the boundary. 
Then, we have 
\begin{eqnarray*}
\omega(f)&=&\frac{1}2(3\sigma(V^4)+e[\xi(\widetilde{f})])\\
&=&\frac{1}2(3\sigma(V^4)+\#\Sigma^{1,1}(\widetilde{f})), 
\end{eqnarray*}
where $\#\Sigma^{1,1}(\widetilde{f})$ denotes 
the algebraic number of cusp points of $\widetilde{f}$. 
\end{thm}

Note that this generalises the result of Hughes and Melvin \cite{h-m} 
since we can consider a usual non-singular Seifert surface for an embedding. 

\subsection{The proofs of Theorems~\ref{thm:proj} and \ref{thm:proj2}}\label{subsect:proof}
\begin{proof}[Proof of Theorem~\ref{thm:proj}] 
In what follows, we often consider $F$ to be an embedding 
in $S^6$ and use the same symbol for it: 
$F\co S^3\hookrightarrow S^6=\Rset^6\cup\{\infty\}$. 

Since the Smale invariant $\omega(f)$ is even, 
we can take the pair $(A,B)$ of integers satisfying 
the following simultaneous linear equations: 

\[\left\{\begin{array}{lcl}
-(A-B)/8&=&\Omega(F)\\
(3A+B)/2&=&\omega(f). 
\end{array}\right.\]

According to Corollary~\ref{cor:signature}, 
take a Seifert surface $\widetilde{F}\co V^4\hookrightarrow S^6$ 
for $F$ with $\sigma(V^4)=A$. 
Then, the Hopf invariant $H_{\widetilde{F}}$ should be equal to $-B$ 
since $\Omega(F)=-(\sigma(V^4)+H_{\widetilde{F}})/8=-(A-B)/8$. 

Consider the normal framing $\nu=(\nu_1,\nu_2,\nu_3)$ 
for $F$ such that the first vector field $\nu_1$ coincides 
with the normal field of $F(S^3)\subset\widetilde{F}(V^4)$. 
This is always possible since $\pi_3(V_{3,1})\approx\pi_3(SO(3))$. 

Then, by using the Compression Theorem \cite{r-s}, 
we can isotope $F$ to a \textit{compressible} embedding $F'$ 
so that the composition $p\circ F'$ is 
an immersion $S^3\looparrowright\Rset^5$. 
Here we mean by ``\textit{compressible}'' that 
the third normal vector field $\nu_3$ 
is parallel to ${\partial}/{\partial{x_6}}$. 
Furthermore, by using a suitable ambient isotopy, 
we can isotope 
$\widetilde{F}\co V^4\hookrightarrow\Rset^6$ to 
$\widetilde{F}'\co V^4\hookrightarrow\Rset^6$ with 
$\partial{\widetilde{F}'}=F'$, 
so that the composition $p\circ\widetilde{F}'\co V^4\to\Rset^5$ has 
no singularity near its boundary. 
Since we can further assume that $p\circ\widetilde{F}'$ is generic, 
$p\circ\widetilde{F}'\co V^4\to\Rset^5$ is a singular Seifert surface 
(see Theorem~\ref{thm:e-s})
for the immersion $p\circ F'\co S^3\looparrowright\Rset^5$. 

Now, we will compute the Smale invariant $\omega(p\circ F')$ 
by using Ekholm-Sz\H{u}cs' formula (Theorem~\ref{thm:e-s}): 
\[
\omega(p\circ F') = \frac{1}{2}(3\sigma(V^4)+e[\xi(p\circ\widetilde{F}')]).
\]
If we denote by $e_{\widetilde{F}'}$ the normal Euler class of 
$\widetilde{F}'$ and by $e(\xi(p\circ\widetilde{F}'))$ 
the Euler class of the bundle $\xi(p\circ\widetilde{F}')$ 
over the set $\widetilde{\Sigma}(p\circ\widetilde{F}')$ 
as in Definition~\ref{defn:euler}, 
then, according to Remark~\ref{rem:cusp} 
(see also \cite[Remark~6.4]{e-s}), 
the homology class of 
$\widetilde{\Sigma}(p\circ\widetilde{F}')$ in $V^4$ 
is dual to $e_{\widetilde{F}'}$ and 
the bundle $\xi(p\circ\widetilde{F}')$ is isomorphic to 
the normal bundle of $\widetilde{F}'$ 
restricted to $\widetilde{\Sigma}(p\circ\widetilde{F}')$. 
Therefore, we see that 
\[
e[\xi(p\circ\widetilde{F}')]
=\langle{e(\xi(p\circ\widetilde{F}')),[\widetilde{\Sigma}(p\circ\widetilde{F}')]}\rangle 
=\langle{e_{\widetilde{F}'}\smile e_{\widetilde{F}'},[V^4,\partial{V^4}]}\rangle, 
\] 
which is further equal to minus the Hopf invariant $-H_{\widetilde{F}}=B$ 
by Theorem~\ref{thm:hopf4k}. 
Finally, we have 
\[
\omega(p\circ F') = \frac{1}{2}(3\sigma(V^4)+e[\xi(p\circ\widetilde{F}')]) 
= \frac{1}{2}(3A+B)
= \omega(f). 
\]
This means that the immersion $f':=p\circ F'$ 
is regularly homotopic to the given immersion $f$. 
\end{proof}

\begin{proof}[Proof of Theorem~\ref{thm:proj2}] 
Take a normal framing $\nu=(\nu_1,\nu_2,\nu_3)$ 
for $F$ such that the third vector field $\nu_3$ is in accordance 
with ${\partial}/{\partial{x_6}}$, which is not tangent to 
$F(S^3)\subset\Rset^6$ since $p\circ F$ is an immersion. 

By Corollary~\ref{cor:normal}, we can consider a Seifert surface 
$\widetilde{F}\co V^4\hookrightarrow S^6$ for which 
the normal field along the boundary 
$F(S^3)\subset\widetilde{F}(V^4)$ 
coincides with the first normal vector field $\nu_1$. 

Then, the composition $p\circ\widetilde{F}$ of $\widetilde{F}$ with the projection $p$ 
has no singularity near the boundary and can be considered 
to be a singular Seifert surface for $p\circ F$. 

By the same argument as in the proof of Theorem~\ref{thm:proj}, 
the Smale invariant $\omega(p\circ F)$ can be computed as: 
\[
\omega(p\circ F) = \frac{1}{2}(3\sigma(V^4)+e_{\widetilde{F}}\smile e_{\widetilde{F}}), 
\]
where $e_{\widetilde{F}}\smile e_{\widetilde{F}}\in H^4(V^4,\partial{V^4})=\Zset$ 
is the square of the normal Euler class of $\widetilde{F}$. 
Since $\sigma(V^4)\equiv e_{\widetilde{F}}\smile e_{\widetilde{F}}\pmod{8}$, 
we have $\omega(p\circ F)\equiv 2\sigma(V^4)\pmod{4}$.
Thus, $\omega(p\circ F)$ is even. 
\end{proof}

\begin{rem}
It seems that Theorem~\ref{thm:proj2} can be deduced 
from the following argument. 
Denote by $\delta(f)$ the number of non-trivial components 
of the set of double points of an immersion 
$f\co S^3\looparrowright\Rset^5$ (see \cite{szucs3}). 
Then, \cite[Proposition~7.5.2 and Corollary~8.3.2]{ekholm1} 
imply that $\delta(f)\pmod{2}$ determines the non-trivial homomorphism 
$\Imm[S^3,\Rset^5]\to\Zset_2$, which consequently coincides with 
the Smale invariant $\omega(f)\pmod{2}$. 
According to Sz\H{u}cs \cite{szucs3}, 
if an immersion $f\co S^3\looparrowright\Rset^5$ 
can be lifted to an embedding in $S^6$, 
then $\delta(f)$ and hence $\omega(f)$ must be even. 
\end{rem}

\section{Punctured embeddings of $4$-manifolds in codimension $2$}\label{sect:punc} 
It is natural to ask if 
a similar construction as in \SubS\ref{subsect:constructions} is 
applicable to other oriented $4$-manifolds. 
We argue here the case of the Kummer surface $K$. 
The following is proved very similarly 
to Proposition~\ref{prop:cp2}. 

\begin{prop}\label{prop:k3}
The punctured Kummer surface $K_\circ$ can be embedded in $S^6$ 
with normal Euler class 
$((0,\ldots,0),(0,\ldots,0),(0,0),(0,0),(2a,2b))\in 
H^2(K_\circ)=H^2(K)=\Zset^8\oplus\Zset^8\oplus\Zset^2\oplus\Zset^2\oplus\Zset^2$ 
\textup{(}represented with respect to the decomposition 
$E_8\oplus E_8\oplus\begin{pmatrix}0&1\\1&0\end{pmatrix}\oplus\begin{pmatrix}0&1\\1&0\end{pmatrix}\oplus\begin{pmatrix}0&1\\1&0\end{pmatrix}$
of the intersection form of $K$\textup{)} for any integers $a,b\in\Zset$. 
Such an embedding restricted 
to the boundary represents $2+ab\in C^3_3=\Zset$. 
\end{prop}

\begin{proof}
Since the Stiefel-Whitney class $w_2(K)$ is 
the reduction of the integral class 
$\eta:=((0,\ldots,0),(0,\ldots,0),(0,0),(0,0),(2a,2b))\in H_2(V^4;\Zset)$, 
by Hirsch's h-principle  \cite{hirsch}, 
there is an immersion of $K_\circ$ in $S^6$ with normal Euler class $\eta$. 
Since the punctured Kummer surface $K_\circ$ has a handlebody decomposition with 
a $0$-handle and several $2$-handles with even framings, 
this immersion is regularly homotopic to an embedding. 

The latter part is just due to a computation 
of the Haefliger invariant: 
\[
-\frac{1}{8}\left(
\sigma(K_\circ)-\begin{pmatrix}2a&2b\end{pmatrix}\begin{pmatrix}0&1\\1&0\end{pmatrix}\begin{pmatrix}2a\\2b\end{pmatrix}
\right)
=\frac{-16-8ab}{-8}
=2+ab. 
\]
This completes the proof. 
\end{proof}

\begin{rem}
In fact, the punctured Kummer surface can be embedded in $S^5$. 
Since such an embedding in $S^5$ has trivial normal bundle, 
its composition with the inclusion $S^5\subset S^6$ restricted 
to the boundary represents $2\in C^3_3=\Zset$ 
(see \cite[Theorem~5.17]{hae2} and \cite[\S7]{takase2}). 
Proposition~\ref{prop:k3}, however, implies that 
for a suitable pair $(a,b)$ of integers, 
the punctured Kummer surface can be a Seifert surface for 
any (isotopy class of) Haefliger knot (see Remark~\ref{rem:any}). 
\end{rem}

\begin{rem}
Combining Corollary~\ref{cor:s2xs2}, Propositions~\ref{prop:cp2}, \ref{prop:k3} 
and Remark~\ref{rem:cp2bar}, 
we see that any indefinite symmetric unimodular form can be 
realised as the intersection form of a Seifert surface 
for a Haefliger knot. 
\end{rem}

We easily see that many other orientable 4-manifolds 
can be Seifert surfaces for Haefliger knots. 
For example, 
the punctured manifold $V^4_\circ$ of 
a closed spin manifold $V^4$ can be embedded in $S^6$, 
since the double $\partial(V^4_\circ\times[0,1])$ of $V^4_\circ$ is 
a closed spin $4$-manifold of signature zero and  
such a manifold can be embedded in $S^6$ 
(e.g., see \cite{c-s,ruberman,saeki} and \cite[Theorem~2.5]{barden}). 
We also see that 
a simply-connected punctured $4$-manifold can be embedded in $S^6$
(\cite[Corollary~4.2]{hirsch2}). 

Note that if $V^4$ is a closed \textit{non-orientable} manifold 
with non-zero third normal Stiefel-Whitney class $\bar{w}_3(V^4)$ 
then its punctured manifold $V^4_\circ$ cannot be embedded in $S^6$. 
For if we have an embedding 
$\widetilde{F}\co V^4_\circ\hookrightarrow S^6$, 
then since $\partial{\widetilde{F}}\co S^3\hookrightarrow S^6$ 
is \textit{topologically} unknotted \cite{stallings}, 
we can construct a \textit{topological} embedding of 
the closed manifold $V^4$ in $S^7$; 
this contradicts $\bar{w}_3(V^4)\ne0$ 
according to \cite[Theorem~1.2]{fang}. 

\ack
The author would like to convey his sincere thanks to 
Professor Osamu Saeki for his invaluable advice.
He would also thank Professors Dennis Roseman, Takashi Nishimura 
and Yukio Matsumoto for constant encouragement 
and many helpful comments. 

The author is partially supported by the Grant-in-Aid for JSPS Fellows.

\end{document}